\documentclass{amsart}
\usepackage{epsfig}
\usepackage{amsfonts}
\usepackage[all,frame,poly]{xy}
\usepackage{graphicx}

\newtheorem{theorem}{Theorem}
\newtheorem{proposition}[theorem]{Proposition}
\newtheorem{lemma}[theorem]{Lemma}
\newtheorem{corollary}[theorem]{Corollary}
\newtheorem{example}{Example}

\def\cal{\mathcal}

\def\BB{\mathbin{\Vert}}

\def\a{\alpha}

\def\bp{\beta_{+}}
\def\bm{\beta_{-}}

\def\gp{\gamma_{+}}
\def\gm{\gamma_{-}}

\def\N{N}
\def\Ntree{N^{\mathrm{tree}}}
\def\Nt#1{N^{{\mathrm{tree}},#1}}
\def\Ntz{\N'}
\def\Ntzprime{\N''}
\def\Nu{N^{\mathrm{uni}}}
\def\Nuu{\Nu_u}
\def\Nb{{\bar N}}
\def\flip#1{\overline{#1}}
\def\floor#1{\lfloor#1\rfloor}
\def\ceil#1{\lceil#1\rceil}
\def\rev#1{\tilde{#1}}
\def\Tree#1{{\cal T}(#1)}

\def\bbb#1{{2#1 \choose #1}}

\date{\today}

\begin{document}
\title{The identity is the most likely exchange shuffle for large $n$}
\author{Daniel Goldstein and David Moews}
\maketitle

\begin{abstract}  Let a deck of $n$ cards be shuffled by successively
exchanging the cards in positions $1$, $2$, \dots, $n$ with cards in randomly
chosen positions. We show that for  $n\ge18,$ the identity permutation 
is the most likely. We prove a 
surprising symmetry of the resulting distribution on permutations.
We also obtain the limiting distribution of the 
number of  fixed points as $n\to\infty.$  
\end{abstract}
\medskip

%\footnotetext{AMS subject classification: 05A15}

\section{Introduction}\label{intro}

One way to shuffle a deck of $n$ cards is as follows. Start by numbering the 
positions of the cards from 1 to $n.$ Swap the first card 
with a card in a randomly chosen position (possibly itself). Then swap
the second card with a card in a randomly chosen position, and
so on; the last exchange swaps the $n$th card  with a card in a
randomly chosen position. Robbins and Bolker \cite{ref4} and Schmidt and Simion
\cite{ref6} have studied the probability distribution on the set of 
permutations of $\{1,\ldots,n\}$ induced by this 
\emph{exchange shuffle}.  
In Section~\ref{goas} we introduce the  directed graph of 
an exchange shuffle. It has vertex set 
$\{1,2,\ldots,n\}$ and one edge for each swap.  A first
result relating the graph to the permutation says that
a connected component of the graph corresponds to 
either a cycle or a product of two disjoint cycles.
In the first case we say that the graph is a tree; 
in the second case we say that the graph is a unicycle because it
contains a unique cycle.
This result is  used in Section~\ref{iimffln} to show
that  for large $n$ the
identity is the most frequent permutation. 
In Section~\ref{roum}, we present a second result
based on a more detailed study of the graph. For a unicyclic graph, this
result says which vertices lie in which cycle, and gives a 
formula for the probability of the permutation.
This  will be used in Sections~\ref{socm} and \ref{fan}. 
In Section~\ref{socm} we
prove various facts about the probability of the concatenation of
a series of permutations. 
In Section~\ref{fan}, we find the most likely permutation in each of the
conjugacy classes of $S_n$.
We prove that 
for $n\ge18$, the identity is most frequent, and
for $4\le n\le17$ the most frequent permutation 
is $(n\ \cdots\ m+1)(m\ \cdots\ 1)$, where $m$ is $n/2$ if $n$ is even,
and either  $(n-1)/2$ or $(n+1)/2$ if $n$ is odd. This answers a question
of Robbins and Bolker.
In Section~\ref{fix} we determine the limiting distribution 
as $n\to\infty$ of the number of fixed points
generated by the exchange  shuffle (and a related shuffle),
thus answering  a question of  Schmidt and Simion.

\section{The graph of a shuffle}\label{goas}

Suppose $\pi$ is a permutation of the set $\{1,\ldots,n\}$.  We wish to count
the number of ways of representing $\pi$ as a product of transpositions
\begin{equation}
\label{eone}
\pi=(n\ a_n) \cdots (2\ a_2) (1\ a_1), \qquad a_1, \ldots, a_n\in\{1,\ldots,n\}.
\end{equation}
Call this number the {\em multiplicity} $\N(\pi)$ of $\pi$.
As the total number of exchange shuffles is $n^n$,
the probability that the permutation $\pi$ arises as an exchange 
shuffle is $N(\pi)/n^n$. In this paper we shall usually consider
multiplicities and not probabilities.

Consider a permutation $\pi$ and its representation (\ref{eone}) as a product 
of transpositions.  We can construct a (possibly loop-containing) directed 
graph (digraph) from this representation by taking $\{1,\ldots,n\}$ 
to be the set of vertices,
and directing an  edge from  $j$ to $a_j$  for each $j$.  We will  
say that this 
digraph also  represents $\pi$.  

Each vertex has out-degree one, so 
every connected component of the graph 
contains a unique directed cycle, possibly with just one
vertex. For, suppose we are given a connected component with vertex set 
of cardinality $m$. Since each vertex has out-degree one, 
the component has exactly $m$
edges; hence there is at most one cycle, directed or otherwise.
Choose a vertex $v_1$  and form a sequence inductively 
by setting $v_{i+1}=a_{v_i}.$ By the pigeon-hole principle this sequence 
repeats.  This repeating part is the unique directed cycle.

To avoid confusion with cycles of permutations, we will call this
directed cycle a {\em ring}. 
The component can also contain trees rooted on this cycle with edges directed
towards the ring.
The transpositions in different 
connected components 
of the digraph do not interact, so our exchange shuffle
induces  a permutation on
the vertex set of each component of the digraph.

\begin{proposition}\label{blah1}
\cite[Lemmas 1.2, 1.4]{ref6}
If the ring has length 1, the permutation $\rho$ induced on
a connected component of the digraph with $m$ vertices 
is an $m$-cycle. 
If the ring has length 2 or greater,
$\rho$  is a product of two disjoint cycles whose lengths sum to $m$.
In the first  case  we call
the component a {\em tree}, in the second  a {\em unicycle}.
\end{proposition}
\begin{proof} 
Suppose the ring has length 1.  Then the component is
a rooted tree with a loop at the root.  The loop has no effect on
the permutation generated. Ignore it, and suppress the corresponding
transposition.  This leaves us with a permutation of the form
\begin{equation}
\label{foo}
\rho = (b_{m-1}\ c_{m-1})(b_{m-2}\ c_{m-2})\cdots (b_1\ c_1), 
\end{equation}
whose  associated graph is a tree  $\cal T$ with vertex set, $\cal S$ say, 
of cardinality  $m$, and
with  edge set $\{(b_i,c_i)| 1\le i\le m-1\}.$

We wish to prove by induction on $m$ that (\ref{foo}) is an $m$-cycle.
For $m=1$ this is true.  For $m$ larger, let $e=(b_i,c_i)$ be any edge
of $\cal T$ with $b_i$  a leaf vertex.  Deleting $e$ and $b_i$
leaves a tree on the vertex set ${\cal S}\setminus\{b_i\}$ of cardinality 
$m-1$, so, by the induction 
hypothesis, the permutation
$$
\rho'=(b_{i-1}\ c_{i-1})
\cdots(b_1\ c_1)(b_m\ c_m)\cdots(b_{i+1}\ c_{i+1})$$
is a cycle.
Write $$\rho'=(c_i\ h_1\ \cdots\ h_{m-2}).$$
Then $\rho$ is conjugate to
$$(b_i\ c_i)\rho'=(c_i\ h_1\ \cdots\ h_{m-2}\ b_i).$$
However, a conjugate of a cycle is a cycle.  This completes the induction.

Next, we treat the case of a ring of length 2 or greater.
Let 
$$
\rho=(f_m\ g_m)\cdots(f_2\ g_2)(f_1\ g_1).
$$
Choose an edge  $(f_i, g_i)$ of the ring.  
Deleting it changes our graph component into a tree, so
$$
\rho'=(f_{i-1}\ g_{i-1})\cdots(f_2\ g_2)(f_1\ g_1)
(f_m\ g_m)\cdots(f_{i+1}\ g_{i+1})
$$
is of the form (\ref{foo}). Then, as we just proved, $\rho'$ is an
$m$-cycle. However, a product of an $m$-cycle and a transposition of two
of its members will equal a product of
two  disjoint cycles whose lengths sum to $m$, thus proving the lemma.
In fact, if we write
$$\rho'=(f_i\ d_1\ \cdots\ d_j\ g_i\ d_{j+1}\ \cdots \ d_{m-2}),$$
$\rho$ will be conjugate to
$$(f_i\ g_i)\rho'=(f_i\ d_1\ \cdots \ d_j)(g_i\ d_{j+1}\ \cdots\ d_{m-2}),$$
a product of two disjoint cycles whose lengths
sum to $m$.  But conjugation does not change cycle structure, so $\rho$
also has this form.
\end{proof}

We can define the multiplicity $\N(\pi)$ of a permutation of any 
set of nonnegative integers
just as we defined the multiplicity of permutations 
of $\{1,\ldots,n\}$.  Write
$\Ntree(\pi)$ for the 
number of representations of $\pi$ as a tree, or {\em tree multiplicity}
of $\pi$, and 
$\Nu(\pi)$ for the number of representations of $\pi$ as a unicycle, or
{\em unicyclic multiplicity} of $\pi$.
Obviously, $\N(\pi)$, $\Ntree(\pi)$, and $\Nu(\pi)$ are invariant under
conjugation of $\pi$ by an order-preserving map.

\begin{lemma}\label{key}
For  a permutation of $\{1,\ldots,n\}$ that is a product 
$\pi=\pi_1\pi_2\cdots\pi_q$ of $q$ disjoint cycles, we have
$$
\N(\pi)=\sum_{\shortstack{\scriptsize $\chi$ an involution of\\
$\scriptstyle \{1,\ldots,q\}$}}\ 
\prod_{\chi(i)=i} \Ntree(\pi_i)
\prod_{\chi(i)=j,\ i<j}  \Nu(\pi_i \pi_j).
$$
\end{lemma}
\begin{proof}Any representation of $\pi$ as an exchange shuffle will
determine a partition of the set $\{\pi_1,\ldots,\pi_q\}$ into
one- and two-element blocks, the one-element blocks being those cycles 
represented as trees, and the
two-element blocks those pairs of cycles represented as unicycles.
Think of such a partition as an involution, by letting
an involution $\chi$ partition
$\{\pi_1,\ldots,\pi_q\}$ into its orbits under $\chi$. 
The product of $\Ntree(\pi_i)$ over the fixed points $\pi_i$ of $\chi$
multiplied by the product of $\Nu(\pi_j\pi_k)$ over the two-element
orbits $\{\pi_j,\pi_k\}$ of $\chi$ will then equal the number of 
representations of $\pi$ yielding this $\chi$, and summing this over
all involutions $\chi$ will yield $\N(\pi)$.
\end{proof}

\section{The identity is most frequent for large $n$}\label{iimffln}
Using Lemma \ref{key} and
the following three lemmas will enable us to prove that the identity 
permutation has highest multiplicity for $n\ge 29$.
In Section~\ref{fan} we will obtain a stronger result.

\begin{lemma} \label{invbd}
Let $Q_n$ be the number of involutions on an $n$-element
set.  Then
\begin{enumerate} 
\item $Q_n 4^{-n}$ is strictly increasing for $n\ge 15$,
\item $Q_n 4^{-n}\le {1\over 4}$ for $n<15$, and 
\item $Q_{29} 4^{-29}>{1\over4}$.
\end{enumerate}
\end{lemma}
\begin{proof}
We have $$Q_0=Q_1=1.$$  For $n\ge 2$, an involution on $\{1,\ldots,n\}$
can either fix the first element, leaving $Q_{n-1}$ ways of moving
the remaining elements, or swap it with a member of $\{2,\ldots,n\}$,
leaving $Q_{n-2}$ ways of moving the remaining elements.  Hence
$$Q_n=Q_{n-1}+(n-1)Q_{n-2}, \qquad \mbox{for\ } n\ge 2.$$
If we set $$R_n=\frac{Q_n}{Q_{n-1}},$$
then induction on $n$ proves that, for all $n\ge 2$,
$$\sqrt{n}<R_n<\sqrt{n}+1.$$
From the above, $Q_n>\sqrt{n} Q_{n-1}$, so
$Q_n 4^{-n}$ will be strictly increasing for $n\ge 15$.  
Inequalities (2) and (3) require a calculation,  which we omit.
\end{proof}

Let $\Nt x(\pi)$ be the number of representations of $\pi$ as a tree
where we constrain $x$ to be the root of the tree, i.e., we constrain
the transposition $(x\ x)$ to be made.  
Let $C_n={1\over n+1}{2n \choose n}$ be the $n$th Catalan number.

% 
% Should we remark about identifying {1, ..., m} with {a_1, ..., a_m}
% as we do here & elsewhere?
%
\begin{lemma} \label{rb} \cite[Theorems 8 and 9]{ref4}  For cycles
$\mu$ on an $n$-element set,
$\Ntree(\mu)=\N(\mu)$ has maximum value $C_n$. Let  $\mu_1$ be the
smallest element of $\mu$; then $\Nt {\mu_1}(\mu)$ has maximum value
$C_{n-1}$. Both these maxima are achieved only
when $\mu$ can be written as $(\mu_n\ \mu_{n-1}\ \cdots\ \mu_1)$,
where $\mu_1<\mu_2<\cdots<\mu_n$.  
\end{lemma}
\begin{proof} 	See Robbins and Bolker.
\end{proof}
\begin{lemma}\label{ouru}
If $\pi$ acts on an $n$-element set, then $\Ntree(\pi)\le 4^{n-1}$
for $n\ge 1$, and $\Nu(\pi)\le 4^{n-2}$ for $n\ge 2$.
\end{lemma}
\begin{proof} We have just seen that $\Ntree(\pi)\le C_n$.  Since $C_1=1$,
the first inequality is true for $n=1$.  We have the following form of
Stirling's formula \cite{ref5}:
$$\sqrt{2\pi m}(m/e)^m e^{1/(12m + 1)} \le m!\le 
\sqrt{2\pi m}(m/e)^m e^{1/12m},\qquad m\ge 1.$$
\medskip
For $n\ge 2$, this  yields the estimate $C_n\le 4^n/(n+1)\sqrt{\pi n}$,
which implies the result.

For the second inequality, let $W_n$ be the maximum of
$\Nu(\pi)$ over all permutations $\pi$ of $\{1,\ldots,n\}$.  We will induce on
$n$ to prove that $W_n\le 4^{n-2}$ for all $n\ge 2$.
For $\Nu(\pi)$ to
be nonzero, $\pi$ must be a product of two disjoint cycles,
$(b_1\ b_2\ \cdots\ b_m)(b_{m+1}\ \cdots\ b_n)$, say.  Without loss of
generality, let $b_m=1$.  Represent $\pi$ in the form (\ref{eone}), and let the
first transposition be $(1\ b_i)$, $i\ne m$.  
The remaining transpositions must represent
$\pi(1\ b_i)$, and they will still do so if we precede them by $(1\ 1)$.
If $i>m$, $\pi(1\ b_i)$ is an $n$-cycle, so the representation is
as a tree, and the tree has a loop at 1.  We thus get at most
$\Nt 1(\pi(1\ b_i))$ choices for the remaining transpositions; by Lemma
\ref{rb}, this is no more than $C_{n-1}$.  If $1\le i\le m-1$, 
$$\pi(1\ b_i)=(b_1\ b_2\cdots\ b_i)(b_{i+1}\ b_{i+2}\cdots\ b_m)
(b_{m+1}\ \cdots\ b_n).$$
Our representation of $\pi(1\ b_i)$ has a loop at 1, so it must contain
a tree component rooted at 1 representing $(b_{i+1}\cdots\ b_m=1)$; for 
the representation of $\pi$ to be unicyclic, it must also contain a unicyclic 
component representing $(b_1\ b_2\cdots \ b_i)(b_{m+1}\ \cdots\ b_n)$, so
by Lemma \ref{rb}, the
number of ways of filling in the remaining transpositions is then no
more than $C_{m-i-1} W_{n-m+i}$.
Then
$$
\Nu(\pi)\le (n-m) C_{n-1} + \sum_{1\le i\le m-1} C_{m-i-1} W_{n-m+i}
$$
and since for all possible $\pi$, $1\le m\le n-1$,
\begin{equation} \label{eee}
W_n \le (n-1) C_{n-1} + \sum_{0\le j\le n-3} C_j W_{n-j-1}.
\end{equation}
$W_n\le 4^{n-2}$
for $2\le n\le 9$ can be checked via direct computer enumeration
(or using some of our results below),
and using the bound (\ref{eee}) recursively 
then yields $W_n\le 4^{n-2}$ for $10\le n\le 21$; so let $n\ge 22$.
Then by induction,
\begin{eqnarray*}
W_n &\le& (n-1) C_{n-1} + \sum_{0\le j\le n-3} C_j 4^{n-j-3}\\
&\le& (n-1) C_{n-1} + 4^{n-3} \sum_{j\ge 0} C_j 4^{-j}.
\end{eqnarray*}
Substitution into the generating function for $C_n$ yields 
$\sum_{j\ge 0} C_j 4^{-j}=2$, and estimation with Stirling's formula
as above yields $(n-1) C_{n-1}\le 2\cdot 4^{n-3}$ for $n-1\ge 64/\pi$,
i.e., $n\ge 22$.  This completes the proof.
\end{proof}

\begin{proposition} 
\label{idbd1}
For $n\ge 29$, the identity has larger multiplicity than any other permutation
of  $\{1,2,\ldots,n\}.$
\end{proposition}
\begin{proof}  
Let $\pi$ be an arbitrary permutation of $\{1,\ldots,n\}$, and let
$\pi=\pi_1\cdots\pi_q$ be its decomposition into disjoint cycles.
By  Lemma~\ref{key},
$$
\N(\pi)=\sum_{\shortstack{\scriptsize $\chi$ an involution of\\
$\scriptstyle \{1,\ldots,q\}$}}\ 
\prod_{\chi(i)=i} \Ntree(\pi_i)
\prod_{\chi(i)=j,\ i<j} \Nu(\pi_i \pi_j).
$$
Plugging in the bounds from Lemma \ref{ouru} above,
$$\N(\pi)\le 
\sum_{\shortstack{\scriptsize $\chi$ an involution of\\
$\scriptstyle \{1,\ldots,q\}$}}\ 4^{n-q},$$ i.e., 
\begin{equation}
\label{ef}
\N(\pi)\le Q_q 4^{n-q}.
\end{equation}
From Lemma \ref{invbd}, it follows that for $n\ge 29$,
$Q_q 4^{-q}$ is maximized at $q=n$, so
$\N(\pi)<Q_n$ unless $q=n$, i.e., unless $\pi$ is the identity permutation.
However the multiplicity of the identity permutation is 
$Q_n$, by \cite[Theorem 6]{ref4}, or  directly from Lemma~\ref{key}.  
This completes the proof.
\end{proof}

\section{A result on unicyclic multiplicity}\label{roum}
The main result of this section, Proposition~\ref{blah2}, is a formula
for the multiplicity of  a product of two disjoint cycles. 

\begin{figure}
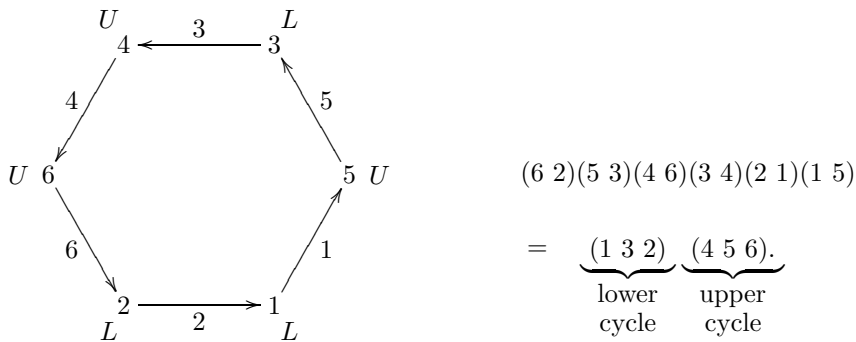

\vbox{

\def\thelab{\ifcase\xypolynode\or 5\or 3\or 4\or 6\or 2\or 1\fi}
\def\thelabb{\ifcase\xypolynode\or U\or L\or U\or U\or L\or L\fi}
\xy
{\xypolygon6"A"{~:{(20,0):}~><{}~*{\thelab}~>>{_{\txt{\thelab}}}}};
{\xypolygon6"B"{~:{(24,0):}~>{}~*{\thelabb}}};
(65,0)*{(6\ 2)(5\ 3)(4\ 6)(3\ 4)(2\ 1)(1\ 5)};
(45,-10)*{=};(57,-10)*+{(1\ 3\ 2)}*\frm{_\}};(71,-10)*+{(4\ 5\ 6).}*\frm{_\}};
(57,-18)*{\txt{lower\\cycle}};(71,-18)*{\txt{\vphantom{lower}upper\\cycle}}
\endxy}
\caption{The lower and upper cycle of a unicyclic component.}
\label{fig1}
\end{figure}

Return to Proposition~\ref{blah1}.
In the case that the connected component is a unicycle, 
the permutation is a product of two disjoint cycles, and we will now
say which ring vertices are in which cycle.  If a ring vertex is larger
than its predecessor, it's in one cycle. If it's smaller, it's in the other.
This is illustrated in Figure~\ref{fig1}.

\begin{proposition}
\cite[Lemma 1.4]{ref6}
The permutation induced on a
digraph component which is a tree is a cycle, and the permutation induced on
a digraph component which is a unicycle is a product of two disjoint
cycles.  For a unicyclic component, divide the vertices in the ring into
{\em upper vertices}, whose  predecessors in the ring have a smaller label, 
and {\em lower vertices}, whose predecessors have a 
larger label.  Then one cycle, 
the {\em upper cycle}, contains all the upper vertices, and the other
cycle, the {\em lower cycle}, contains all the lower vertices.
\end{proposition}
\noindent{{\it Remark 1.}  We will find it useful from now on to
label edges in our digraphs by their source vertices.  This was already
seen in Figure 1.}\hfill\\
{{\it Remark 2.}  Later, in Proposition~\ref{blah2}, we will see which 
non-ring vertices are upper and which are lower. 
We will also see what order to take them in.}\hfill\\ 
\begin{proof} 
First, let the component be a tree. We will prove the induced permutation
is a cycle, but now without using conjugation.  As before, let us suppose we 
have a permutation of the form (\ref{foob}):
\begin{eqnarray}
\label{foob}
\rho &=& (b_{m-1}\ c_{m-1})(b_{m-2}\ c_{m-2})\cdots (b_1\ c_1), \\
\qquad && 
       \hbox{whose associated graph, $\cal T,$ is a tree with vertex set 
$\cal S$ of}
\nonumber \\
\qquad && \hbox{cardinality $m$ 
and edge  set $\{(b_1,c_1),\ldots,(b_{m-1},c_{m-1})\}$.}\nonumber
\end{eqnarray}

We wish to prove by induction on $m$ that (\ref{foob}) is an $m$-cycle.
We will need a stronger induction hypothesis.
Let $b_1$, \dots, $b_{m-1}$ be thought of as an increasing sequence of times 
of day, and let a postman travel on the set $\cal S$.  If he is at $b_i$
at time $b_i$, he will move to $c_i$, or, if he is at $c_i$ at time $b_i$,
he will move to $b_i$.  Otherwise, he will stay fixed.  Our induction
hypothesis is that, regardless of where and at what time of day the postman 
starts, he will return to the same place at the same time of day $m$ days later,
and that his intervening journey will, for any time of day $t$ and 
$s\in \cal S$, find him at $s$ at time $t$ exactly once.  

For $m=1$ this is clear,
since the postman simply stays put for one day.  Otherwise, let $e=(b_i,c_i)$
be any edge of $\cal T$.  Deleting $e$ will split $\cal T$ into subtrees;
let the subtree where the postman starts his journey be $\cal T'$,
with $m'$ vertices, say, and let the other subtree be $\cal T''$, 
with $m''=m-m'$ vertices.  
By the induction hypothesis on $\cal T'$, the postman will travel until he 
reaches an endpoint of $e$ at time of day $b_i$, when he will travel to the 
other endpoint of $e$ in $\cal T''$.  Then,
by the induction hypothesis on $\cal T''$, he will find himself at the endpoint
of $e$ at time of day $b_i$ exactly $m''$ days later,
when he will travel back to 
$\cal T'$.  His journey there will then go on as if he had never left, so
he will arrive back at his starting point and starting time of day
after having spent $m'$ days in $\cal T'$ and $m''$ days in $\cal T''$, for a 
total of $m'+m''=m$ days.  In fact, his cyclical journey 
in $\cal T$ simply consists of his journeys in $\cal T'$ and $\cal T''$ 
spliced together, so it's clear that the induction hypothesis is still 
satisfied.

If we let the postman be at $v$ at the start of the day (before any time $b_i$),
then he will be at $\rho(v)$ at the end of the day (after all times $b_i$.)
We may say that he stops for the night there.  At any rate, he is still at 
$\rho(v)$ at the start of the next day;  he is at $\rho^2(v)$ at the 
start of the day after that; and so on.  
Our induction hypothesis then tells us that $v$, 
$\rho(v)$, \dots, $\rho^{m-1}(v)$ are distinct, so $\rho$ is an $m$-cycle, as 
desired.

If the component is unicyclic, suppose at first that there are no edges 
other than those in the ring, and draw the ring in the plane so that it 
is oriented counterclockwise.  If the postman spends the night at an upper 
vertex, he will leave via the incoming edge, in a clockwise direction,
and continue to travel clockwise, passing through all the lower 
vertices he meets on the way, before spending the next night at the next 
upper vertex.  Therefore, the upper vertices form a cycle in the induced
permutation.  Similarly, if he spends a night at a lower vertex, the postman
will spend his next night at the next lower vertex in a counterclockwise 
direction, so the lower vertices also form a cycle.  We illustrate this in 
Figure 1.
In the general case, the unicyclic component will contain trees attached 
to the ring by edges.  However, we already know that each such tree will, on 
its own, generate a cyclical journey.  When attached to the ring, these 
journeys will be spliced into the journeys around the cycles of upper and lower vertices, as explained above in the proof that the induced permutation on
a tree is a cycle.  This completes the proof.
\end{proof}

We wish  to state a formula for $\Nu$.
Let $\pi$ be a permutation on a set of nonnegative integers that contains 
zero.
Write $\Ntz(\pi)$ for $\Nt 0(\pi)$, the number of representations of $\pi$ as 
a  tree with root 0.  
Given a tree $T$ with
root 0, let $\bp(T)$ be the largest vertex adjacent to 0, and $\bm(T)$ the
smallest.  Then we write $\Ntz_{\bp>x}(\pi)$ for the number of representations
of $\pi$ as a tree $T$ with root 0 and $\bp(T)>x$;
$\Ntz_{\bm<x}(\pi)$ is defined similarly.

We will use capital letters $A$, $B$, $D$, and $E$ to denote nonempty
sequences of distinct nonnegative integers; the concatenation of 
sequences $A$ and $B$ will
be written $A\ B$.  
Let $\gm(A)$ be the first element of $A$, and $\gp(A)$ the
last.  Let $\Nuu((A),(B))$ be the number of representations
of the permutation $(A)(B)$ as a unicyclic digraph for which
$(A)$ is the upper cycle.  We can now state our formula for $\Nu$.
\begin{proposition}\label{blah2}
\label{claim1}
For all disjoint sequences of distinct positive integers $A$ and $B$,
$\Nuu((A),(B))$ is equal to 
\begin{equation}\label{eqm}
\sum_{k\ge 1}
\frac{1}{k}
\sum_{\shortstack{$\scriptstyle (A_k\cdots A_1) = (A)$\\
$\scriptstyle (B_1\cdots B_k)= (B)$}}
\prod_i \Ntz_{\bp>\gm(B_{i})}((A_i\ 0))
\prod_j \Ntz_{\bm<\gp(A_{j+1})}((B_j\ 0)),
\end{equation}
where we take $A_{k+1}=A_1$.
\end{proposition}
\begin{proof}
First we  comment on what is to be proved.
Start with the left-hand side above.
Suppose we are given a representation of $(A)(B)$ as a unicyclic
exchange shuffle, 
and suppose that in the associated directed graph $(A)$ is the upper and
$(B)$ the lower cycle. 
Break up the ring into maximal 
segments of consecutive vertices
that are all upper or all lower. Say there are $2k$ such segments. 
Beginning with an
arbitrary upper segment, write the vertices of the ring as 
 
$$\mathrm{Ring } = D_1\  E_1\ D_2\ E_2\ \cdots\  D_k\  E_k.$$
Of course, there are $k$ ways of doing this, depending on which upper segment 
we begin with.

We construct a forest of $k$ trees by, for each $i=1$, \dots, $k$, doing
the following:
\begin{enumerate}
\item Start with the digraph representing $(A)(B).$
\item Follow the postman as he traverses the part $D_i$ of the ring,
and certain trees hanging off the ring (see Figure~2).
\item Use the dummy vertex 0 for any ring vertices the postman traverses 
that are not in $D_i.$
\end{enumerate}

In this way, we get upper trees $T_i$ for $1\le i \le k$. 
By a similar process, to be described below,
we get lower trees $U_j$ for $1\le j \le k$. 
We get  corresponding decompositions
$$A = A_k\ A_{k-1}\ \cdots\ A_1$$
and
$$B = B_1\ B_2\ \cdots\ B_k.$$
The details of these decompositions are illustrated in Figures~2 and 3.
These $2k$ trees and $2k$ sequences
satisfy the properties implied by the right-hand side above:

\begin{enumerate}
\item The tree  $T_i$ represents $(A_i\ 0).$ 
The tree  $U_j$ represents $(B_j\ 0).$
\item The largest neighbor of $0$ in the tree
$T_i$  is greater than the first element of $B_i$.
The smallest  neighbor of $0$ in the tree
$U_j$  is less  than the last element of $A_{j+1}$.
\end{enumerate}
We let these indices wrap around modulo $k$, so that $A_{k+1}=A_1$.
 
The content of Theorem~\ref{blah2} is that this decomposition exists, 
and, conversely, that any choice of
positive integer $k$ and choice of $2k$ subsequences and
$2k$ tree representations thereof  satisfying  conditions 1 and 2 above 
give rise to a unicyclic representation of $(A)(B)$ where $(A)$ is the
upper and $(B)$ the lower cycle.

Next we proceed to the proof. Consider a representation of $(A)(B)$ as
a unicycle, where $(A)$ is the upper cycle.  Take the vertices of
the ring to be, in order, $\a_1$, $\a_2$, \dots, 
$\a_m$, where $\a_1$ is upper and $\a_m$ lower.
We will let these indices wrap around modulo $m$, so $\a_0=\a_m$, 
$\a_{m+1}=\a_1$, and so on.  Start with a maximal consecutive subsequence 
$D_i= \a_a\ \a_{a+1}\ \cdots\ \a_{a+b}$ 
of upper vertices among the $\a_i$'s, and let
$\a_{a+b+1}\ \cdots\ \a_{a+b+c}$ be the following maximal consecutive
subsequence of lower vertices.
Suppose that our postman has just started out for the day from 
$\a_{a+b+c+1},$ which is an upper vertex, by maximality.
What happens to him in the remainder of his journey?

Since $\a_{a+b+c+1}$ is upper, $\a_{a+b+c+1}>\a_{a+b+c}$, so he
will eventually be leaving along the edge labeled $\a_{a+b+c}$; but first,
he may traverse subtrees of the unicyclic component attached to the ring at
$\a_{a+b+c+1}$.
Call a collection of subtrees attached to the ring by edges labeled $j$,
where $\kappa<j<\lambda$, a $\Tree{\kappa,\lambda}$.
We can now say that our postman will 
traverse the $\Tree{-\infty, \a_{a+b+c}}$ incident on $\a_{a+b+c+1}$.  
Following this, he arrives at $\a_{a+b+c}$.  This is a lower vertex,
so $\a_{a+b+c-1}>\a_{a+b+c}$.  He will therefore traverse the incident
$\Tree{\a_{a+b+c},\a_{a+b+c-1}}$ and continue to $\a_{a+b+c-1}$
without stopping for the night.  If this vertex is a lower vertex, he
will traverse $\Tree{\a_{a+b+c-1},\a_{a+b+c-2}}$ and continue to 
$\a_{a+b+c-2}$, and so on.  Eventually, after traversing
$\Tree{\a_{a+b+1}, \a_{a+b}}$, he will come to $\a_{a+b}$.  This
is an upper vertex, so $\a_{a+b-1}<\a_{a+b}$, meaning that he will
traverse the incident $\Tree{\a_{a+b},\infty}$, and then stop for 
the night at $\a_{a+b}$.  The next day, he will traverse the 
$\Tree{-\infty, \a_{a+b-1}}$ before traveling to $\a_{a+b-1}$.
Here he will traverse $\Tree{\a_{a+b-1},\infty}$, stop at $\a_{a+b-1}$,
and traverse $\Tree{-\infty,\a_{a+b-2}}$ before passing to $\a_{a+b-2}$,
and he will continue in this way until reaching $\a_a$ and stopping for
the night.  Since we are now in the same position we started out in
(although with different values of $a$, $b$, and $c$), we end our
description here.
This travelogue generates a sequence which forms a portion of the cycle 
$(A)$.  Equating our $\Tree{\kappa,\lambda}$'s 
with the sequence of vertices they generate when traversed, we may write it as
\begin{eqnarray}
\label{eqai}
A_i&=&
\Tree{-\infty, \a_{a+b+c}} \qquad \mbox{(rooted at $\a_{a+b+c+1}$)}\\
&&\Tree{\a_{a+b+c},\a_{a+b+c-1}} \qquad \mbox{(rooted at $\a_{a+b+c}$)}
\nonumber \\
&&\vdots\nonumber \\
&&\Tree{\a_{a+b+1}, \a_{a+b}}
\qquad \mbox{(rooted at $\a_{a+b+1}$)}
\nonumber \\
&&
\Tree{\a_{a+b},\infty}\  \a_{a+b}\  \Tree{-\infty, \a_{a+b-1}}
\qquad \mbox{(both trees rooted at $\a_{a+b}$)}
\nonumber \\
&&\vdots \nonumber \\
&&
\Tree{\a_{a+1},\infty}\  \a_{a+1}\  \Tree{-\infty, \a_a}
\qquad \mbox{(both trees rooted at $\a_{a+1}$)}
\nonumber \\
&&\Tree{\a_a, \infty}\  \a_a
\qquad \mbox{(tree rooted at $\a_a$)}. \nonumber
\end{eqnarray}

Generate a tree $T_i$ by the following procedure: take the vertices 
$\a_a$, $\a_{a+1}$, \dots, $\a_{a+b}$, and the edges connecting them;
add to this all the $\Tree{\kappa,\lambda}$'s above rooted on these 
vertices.  Add an edge from $\a_{a+b}$ to 0, and take all the trees
$\Tree{-\infty, \a_{a+b+c}}$,
$\Tree{\a_{a+b+c},\a_{a+b+c-1}}$,
\dots,
$\Tree{\a_{a+b+1}, \a_{a+b}}$ appearing above and reroot them at 0.
Walking the postman through $T_i$ tells us immediately that
it is a representation of $(A_i\ 0)$.  Observe that $\bp(T_i)=\a_{a+b}$,
and that $\gp(A_i)=\a_a$.
All this is illustrated in Figure 2.

Let us interchange `lower' and `upper' above, so that
our starting subsequence is a maximal consecutive subsequence 
$E_j=\a_d\ \a_{d+1}\ \cdots\ \a_{d+e}$ of lower vertices,
$\a_{d+e+1}$, \dots, $\a_{d+e+f}$ are upper, and $\a_{d+e+f+1}$ is lower.
Supposing that our postman is just about to stop for the night at
$\a_d$, he will travel along the sequence
\begin{eqnarray}
\label{eqbj}
B_j&=&
\a_d\  \Tree{-\infty, \a_d}\qquad \mbox{(tree rooted at $\a_d$)}\\
&&\Tree{\a_d,\infty}\  \a_{d+1}\  \Tree{-\infty, \a_{d+1}}
\qquad \mbox{(both trees rooted at $\a_{d+1}$)}\nonumber \\
&&\vdots\nonumber\\
&&\Tree{\a_{d+e-1},\infty}\  \a_{d+e}\  \Tree{-\infty, \a_{d+e}}
\qquad \mbox{(both trees rooted at $\a_{d+e}$)}\nonumber \\
&&\Tree{\a_{d+e},\a_{d+e+1}} \qquad \mbox{(rooted at $\a_{d+e+1}$)}
\nonumber \\
&&\vdots\nonumber \\
&&\Tree{\a_{d+e+f-1},\a_{d+e+f}} \qquad \mbox{(rooted at $\a_{d+e+f}$)}
\nonumber \\
&&\Tree{\a_{d+e+f},\infty} \qquad \mbox{(rooted at $\a_{d+e+f+1}$)}\nonumber
\end{eqnarray}
before stopping for the night at $\a_{d+e+f+1}$, where we will leave him.
As before, we may view $B_j$ as being part of $(B)$. 
Generate a tree
$U_j$ by taking the vertices $\a_d$, $\a_{d+1}$, \dots, $\a_{d+e}$,
and the edges connecting them; add to this all the $\Tree{\kappa,\lambda}$'s 
above rooted on these vertices.  
Add an edge from $\a_{d+e}$ to 0, and take all the trees
$\Tree{\a_{d+e}, \a_{d+e+1}}$,
\dots,
$\Tree{\a_{d+e+f-1}, \a_{d+e+f}}$,
$\Tree{\a_{d+e+f},\infty}$ appearing above and reroot them at 0.
Walking the postman through $U_j$ tells us that
it is a representation of $(B_j\ 0)$.  In this case, $\bm(U_j)=\a_{d+e}$,
and $\gm(B_j)=\a_d$. 
This is illustrated in Figure 3.

Split the vertices occurring in the ring into blocks
of upper and lower vertices, so that
$$\a_1\ \cdots\ \a_m=D_1\ E_1\ D_2\ E_2\ \cdots\ D_k\ E_k$$
for some sequences $D_1$, \dots, $D_k$ of upper and $E_1$, \dots,
$E_k$ of lower vertices.  We can then apply the above transformations
to generate $A_1$, \dots, $A_k$, $B_1$, \dots, $B_k$ and associated
tree representations $T_i$ of the $(A_i\ 0)$'s and $U_j$ of the
$(B_j\ 0)$'s.  Looking at the postman's journey tells us that
$(A)=(A_k\ \cdots\ A_1)$ and 
$(B)=(B_1\ \cdots\ B_k)$.
As above, let 
the first vertex of $D_i$ be $\a_a=\gp(A_i)$.  This must be upper, so 
$\a_{a-1}<\a_a$; but $\a_{a-1}$ is the last vertex of $E_{i-1}$,
which is $\bm(U_{i-1})$ (letting indices wrap around modulo $k$),
so $\bm(U_{i-1})<\gp(A_i)$.  Similarly, $\bp(T_{i})>\gm(B_i)$.

In this way, we get a split of $(A)$ and $(B)$ into subsequences,
and tree representations thereof, of exactly the kind counted on
the right-hand side of (\ref{eqm}).  In fact, for each unicyclic 
representation whose ring contains $k$ blocks of upper vertices interleaved 
with $k$ blocks of lower vertices, we get $k$ such splits and representations, 
since we had $k$ choices for $\a_1$.  This gives a map from the set of
things counted on the left-hand side of (\ref{eqm}) to the set of 
things counted on the right-hand side of (\ref{eqm}).

To find an inverse map, we need to be able to splice a number of 
tree representations together into a unicycle.  Suppose that we start
with a tree $T_i$ which represents $(A_i\ 0)$.  Our first task is to
construct the corresponding $D_i=\a_a\ \cdots\  \a_{a+b}$, say.
We may let $\a_{a+b}$ be the largest vertex with an edge to 0.  
Set $\a_{a+b-1}$
to be  the largest vertex less than $\a_{a+b}$ with an
edge to $\a_{a+b}$, set $\a_{a+b-2}$ to the largest vertex less than 
$\a_{a+b-1}$
with an edge to $\a_{a+b-1}$, and so on, until we reach a vertex $\a_a$ such 
that no $x<\a_a$ has an edge to $\a_a$.  The sequence of vertices constructed
in this way obviously has $\a_a<\a_{a+1}<\cdots<\a_{a+b}$. The tree
$T_i$ will consist of the path of edges from $\a_a$ to 0, together with
a $\Tree{\a_a,\infty}$ rooted at $\a_a$,
a $\Tree{-\infty,\a_a}$ and a $\Tree{\a_{a+1},\infty}$ rooted at 
$\a_{a+1}$,
a $\Tree{-\infty,\a_{a+1}}$ and a $\Tree{\a_{a+2},\infty}$ rooted at 
$\a_{a+2}$, and so on, until we reach 0, which will have a
$\Tree{-\infty,\a_{a+b}}$ rooted at it.  
We therefore have $\gp(A_i)=\a_a$ and $\bp(T_i)=\a_{a+b}$.

Similarly, if we start with a tree $U_j$ representing $(B_j\ 0)$,
we need to construct $E_j=\a_d\ \cdots\ \a_{d+e}$.
Let $\a_{d+e}$ be the smallest vertex with an edge to 0, $\a_{d+e-1}$
the smallest vertex greater than $\a_{d+e}$ with an
edge to $\a_{d+e}$, $\a_{d+e-2}$ the smallest vertex greater than $\a_{d+e-1}$
with an edge to $\a_{d+e-1}$, and so on, until we reach a vertex $\a_d$ such 
that no $x>\a_d$ has an edge to $\a_d$.  In this case,
$\a_d>\a_{d+1}>\cdots>\a_{d+e}$, $\gm(B_j)=\a_d$, and $\bm(U_j)=\a_{d+e}$.

Now that we have constructed the $D_i$'s and $E_j$'s, we can concatenate
them into a ring by placing an edge from the last vertex of each $D_i$ to
the first vertex of $E_i$, and the last vertex of $E_i$ to the first of
$D_{i+1}$.  We delete the edges from the last vertices of the $D_i$'s and
$E_i$'s to 0, and leave the trees incident upon the vertices of the
$D_i$'s and $E_i$'s in place.

Let $\a_a$ be the first vertex of $D_i$.  The inequality 
$\gp(A_i)>\bm(U_{i-1})$ then tells us that $\a_a$ exceeds the last vertex
of $E_{i-1}$, i.e., $\a_a>\a_{a-1}$.  Similarly, if $\a_d$ is the first
vertex of $E_i$, the inequality $\gm(B_i)<\bp(T_i)$ tells us that
$\a_d<\a_{d-1}$.  This tells us that the vertices we wish to be upper
are indeed upper, and the ones we wish to be lower, lower.
The $\Tree{-\infty, \a_{a+b}}$ of $T_i$ rooted at 0
must be unrooted and rerooted on the vertices of $E_i$ and the first vertex
of $D_{i+1}$, so as to represent $A_i$ in the form (\ref{eqai}); 
similarly, the $\Tree{\a_{d+e}, \infty}$
of $U_j$ rooted at 0 must be rerooted on the vertices of $D_{j+1}$ and
the first vertex of $E_{j+1}$, so that $B_j$ is represented in the
form (\ref{eqbj}).  From the inequalities we have established, we see
that this can be done in a unique way.  This construction of a
unicyclic representation of $(A)(B)$ shows us that we have a map from the
set of things counted by the right-hand side of (\ref{eqm}) to the set of
things counted by its left-hand side.  This map is inverse to our previous
map, so both maps are bijections, and the left- and right-%
hand sides of (\ref{eqm}) are equal.  This completes the proof.
\end{proof}

\begin{corollary} If every member of $A$ is less than every member of $B$,
then
$$\Nuu((A),(B))=0$$ and
\begin{equation}
\label{eqo}
\Nu((A)(B))
=\Nuu((B),(A))=
\sum_{k\ge 1} 
\frac{1}{k}
\Nb_k((A)) \Nb_k((B)),
\end{equation}
%so
%\begin{equation}
%\label{eqo}
%\Nu((A)(B))=\Nuu((A),(B))+\Nuu((B),(A))=
%\sum_{k\ge 1} 
%\frac{1}{k}
%\Nb_k((A)) \Nb_k((B)),
%\end{equation}
where
$$
\Nb_k((D))=\sum_{(D_1\cdots D_k)= (D)} \prod_i \Ntz((D_i\  0)).
$$
\end{corollary}
\begin{proof}
If every member of $A$ is less than every member of $B$, then
$\Ntz_{\bp>\gm(B_{i})}((A_i\ 0))$ must always vanish, yielding 
$\Nuu((A),(B))=0$.  On the other hand, if
$T$ is a tree with root 0 
representing $(B_i\ 0)$, then the condition $\bp(T)>\gm(A_i)$ is vacuous,
so $$\Ntz_{\bp>\gm(A_i)}((B_i\  0))=\Ntz((B_i\  0)).$$  Similarly,
$$\Ntz_{\bm<\gp(B_{j+1})}((A_j\  0))=\Ntz((A_j\  0)),$$ and substituting
into (\ref{eqm}) yields the desired result.
Of course, the first equation is obvious and does not require
 Proposition~\ref{blah2}.
\end{proof}

\section{Symmetries of concatenation multiplicities}\label{socm}

If $\pi$ is a permutation on $\{1,\ldots,n\}$ and
$\rho$ is a permutation on $\{1,\ldots,m\}$, let the 
{\em concatenation}
of $\pi$ and $\rho$, $\pi\BB\rho$,
be the permutation on $\{1,\ldots,m+n\}$ that 
satisfies $(\pi\BB\rho)(j)=\pi(j)$ for $j=1$, \dots, $n$, and
$(\pi\BB\rho)(j)=n+\rho(j-n)$ for $j=n+1$, \dots, $m+n$.
Also, let $\sigma_n$ be the (self-inverse) permutation 
on $\{1,\ldots,n\}$ sending 
each $i$ to $n+1-i$, and let $\flip{\pi}$, the {\em flip} of $\pi$, be 
$\sigma_n \pi^{-1} \sigma_n$.
Since we have
$$
\flip{ (n\ a_n) \cdots (2\ a_2) (1\ a_1)}=
(n\  \sigma_n(a_1)) \cdots (2\  \sigma_n(a_{n-1})) (1\  \sigma_n(a_n)),
$$
there is a natural bijection between the representations of $\pi$ and
those of $\flip{\pi}$; therefore, 
\begin{equation}
\label{eq1}
\N(\pi)=\N(\flip{\pi}).
\end{equation}

Our formula for $\Nu$ now yields some symmetry properties of 
the multiplicity of a concatenation.

\begin{theorem}
\label{resa}
$\N(\rho_1\BB\rho_2\BB\rho_3\BB\rho_4)=
\N(\rho_1\BB\rho_3\BB\rho_2\BB\rho_4)$, 
      for all (possibly empty) $\rho_1$, $\rho_2$, $\rho_3$, and $\rho_4$. 
\end{theorem}

\begin{theorem}
\label{resb}
$N(\flip{\rho_1}\BB\rho_2)=\N(\rho_1\BB\rho_2)$, for all 
$\rho_1$ and $\rho_2$.
\end{theorem}

It is a corollary of these results that, in a concatenation 
$\pi_1\BB\cdots\BB\pi_q$
, the permutations making up the concatenation can have
their order rearranged, or be flipped, without affecting the 
overall multiplicity.
\smallskip

\noindent{\it Proof of Theorem \ref{resa}}. To 
prove Theorem \ref{resa}, let $\pi=(\rho_1\BB\rho_2\BB\rho_3\BB\rho_4)$
be decomposed into disjoint cycles, $\pi=\pi_1\cdots\pi_q$.  Then
there is a corresponding decomposition 
$\pi'=(\rho_1\BB\rho_3\BB\rho_2\BB\rho_4)=\pi'_1\cdots\pi'_q$,
where each $\pi'_i$ is a conjugate of $\pi_i$ by an additive translation.
Now since, by Lemma~\ref{key},
$$
\N(\pi)=\sum_{\shortstack{\scriptsize $\chi$ an involution of\\
$\scriptstyle \{1,\ldots,q\}$}}\ 
\prod_{\chi(i)=i} \Ntree(\pi_i)
\prod_{\chi(i)=j,\ i<j} \Nu(\pi_i \pi_j).
$$
and similarly for $\pi'$, it will do to show that 
\begin{equation}
\label{glorb1}
\Ntree(\pi_i)=\Ntree(\pi'_i) \qquad\hbox{ for all $i$}
\end{equation}
and 
\begin{equation}
\label{glorb2}
\Nu(\pi_i \pi_j)=\Nu(\pi'_i\pi'_j)\qquad\hbox{ for all distinct $i$ and $j$.}
\end{equation}
However, $\pi'_i$ and $\pi_i$ are conjugate by an additive translation, 
and additive 
translation is order-preserving.  This implies (\ref{glorb1}).
As for (\ref{glorb2}), if 
$\pi_i$ and $\pi_j$ come from the same $\rho_k$, $\pi'_i \pi'_j$ and
$\pi_i \pi_j$ are conjugate by an additive translation, so the equality 
follows.  Otherwise, every element in the domain of $\pi_i$ is smaller than
every element in the domain of $\pi_j$, or vice-versa, so by (\ref{eqo}), 
$$\Nu(\pi_i\pi_j)=\sum_{k\ge 1}{1\over k} 
\Nb_k(\pi_i) \Nb_k(\pi_j),$$
and similarly for $\Nu(\pi'_i\pi'_j)$.  Conjugacy by an additive
translation then 
implies that $\Nb_k(\pi_i)=\Nb_k(\pi'_i)$ for all $i$, proving 
(\ref{glorb2}) and Theorem~\ref{resa}.
\qed\medskip

Here is an alternate proof of (\ref{glorb2}) that does not use
(\ref{eqo}). First, look at a concatenation of two
permutations, $\pi$  a permutation on $\{1,\dots,n\}$ and 
$\rho$ a permutation on $\{1,\dots,m\}$.
In Lemma~\ref{lem12}, below, we
give an explicit bijection between the sets of representations of
$(\pi\BB\rho)$ and $(\rho\BB\pi)$.  If $\rho$ and $\pi$ are cycles, then,
by Lemmas \ref{key} and \ref{lem12},
\begin{eqnarray*}
\Ntree(\pi)\Ntree(\rho)+\Nu(\pi\BB\rho)
&=&
\N(\pi\BB\rho)
\\
&=&\N(\rho\BB\pi)\\
&=&\Ntree(\rho)\Ntree(\pi)+\Nu(\rho\BB\pi),
\end{eqnarray*}
so $\Nu(\pi\BB\rho)=\Nu(\rho\BB\pi)$.  Together with the independence
of $\Nu$ under conjugation by an order-preserving map, this implies 
(\ref{glorb2}).
\begin{lemma}
\label{lem12}
The two permutations 
$(\pi\BB\rho)$ and $(\rho\BB\pi)$ of $\{1,\dots,m+n\}$ have the same 
multiplicity $\N(\pi\BB\rho) = \N(\rho\BB\pi)$. Moreover, this equality comes 
from an explicit bijection of representations.
\end{lemma}
\begin{proof}
Represent  $(\pi\BB\rho)$ as a product of transpositions,
\begin{equation}
(\pi\BB\rho) = (m+n\ a_{m+n})\cdots(1\ a_1).\label{ep1}
\end{equation}
We will exhibit a sequence $(b_1,\ldots,b_{m+n})$ with 
$$
(\rho\BB\pi) = (m+n\ b_{m+n})\cdots(1\ b_1).
$$
Let $$\psi(x)=\left\{
\begin{array}{ll}
m+x & \mbox{if  } 1\le x\le n\\
x-n & \mbox{if  } n+1\le x\le m+n
\end{array}
\right. 
$$
be the permutation of  $\{1,\dots,m+n\}$ that moves the
first  $n$ terms past the  last $m$. Then by definition 
$$(\pi\BB\rho) = \pi\cdot \psi^{-1}\rho\psi,$$ and similarly
$$(\rho\BB\pi) = \rho\cdot \psi\pi\psi^{-1}.  $$
Using $\psi$, (\ref{ep1}) becomes
$$
\pi\cdot \psi^{-1}\rho\psi = [(m+n\ a_{m+n})\cdots(n+1\ a_{n+1})]
[(n\ a_{n})\cdots(1\ a_{1})]
$$
Rearrange terms to get:
$$
\mbox{id} = [\psi\pi^{-1}(m+n\ a_{m+n})\cdots(n+1\ a_{n+1})]
[(n\ a_{n})\cdots(1\ a_{1})\psi^{-1}\rho^{-1}].
$$
Now we use repeatedly the conjugation rule 
$\tau (x\ y)\tau^{-1} = (\tau(x)\ \tau(y))$ to obtain 
$$
\mbox{id} = [(m\ b_m)\cdots(1\ b_1)\psi\pi^{-1}] 
[\psi^{-1}\rho^{-1}(m+n\ b_{m+n})\cdots(m+1\ b_{m+1})]
$$
where 
\begin{equation}
\label{themap}
b_x=\left\{
\begin{array}{ll}
\psi\pi^{-1}(a_{x+n}) & \mbox{if  } 1\le x\le m\\
\rho\psi(a_{x-m}) & \mbox{if  }m+1\le x\le m+n.
\end{array}
\right. 
\end{equation}
Now rearrange again to get 
$$(\rho\BB\pi) = \rho\cdot\psi\pi\psi^{-1}= (m+n\ b_{m+n})\cdots(1\ b_1).$$
This shows us that (\ref{themap})
induces a map from the representations of $(\pi\BB\rho)$ to those of
$(\rho\BB\pi)$.  Interchanging the roles of $m$
and $n$, and of $\rho$ and $\pi$, in (\ref{themap}) will give us a map from
representations of $(\rho\BB\pi)$ to those of $(\pi\BB\rho)$, and it
will be inverse to this map.  This shows that the pairing
$$
(a_1,\ldots,a_{m+n})\longleftrightarrow(b_1,\ldots,b_{m+n})
$$
yields a bijection of representations of $(\pi\BB\rho)$ and $(\rho\BB\pi)$.
\end{proof}

\noindent{\it Proof of Theorem \ref{resb}}. Let $\rho_1$ act on 
$\{1,\ldots,n\}$ and $\rho_2$ act on $\{1,\ldots,m\}$.
If we go through the same analysis as in Theorem \ref{resa}, 
we will see that it 
suffices to prove that, for any
disjoint cycles $\pi$ and $\mu$ with domain contained in $\{1,\ldots,n\}$, 
and any cycle $\nu$ with domain contained in $\{n+1,\ldots, m+n\}$, that
\begin{equation}
\label{e1}
\Ntree(\sigma_n \pi^{-1} \sigma_n)=\Ntree(\pi),
\end{equation}
\begin{equation}
\label{e2}
\Nu(\sigma_n (\pi\mu)^{-1} \sigma_n)=\Nu(\pi\mu),
\end{equation}
and
\begin{equation}
\label{e3}
\Nu(\sigma_n \pi^{-1} \sigma_n \nu)=\Nu(\pi\nu).
\end{equation}

Let a permutation $\phi$ have domain $\{e_1<e_2<\cdots<e_r\}\subseteq
\{1,\ldots,n\}$.  Now $\sigma_n$ is order-reversing on the domain of $\phi$,
so there is a natural bijection between the representations
of $\phi$ and those of $\sigma_n \phi^{-1} \sigma_n$, given by sending
the representation
$$
\phi = (e_r\  e_{a_r})\cdots(e_2\ e_{a_2})(e_1\ e_{a_1})
$$
to the representation
$$
\sigma_n \phi^{-1} \sigma_n = 
(\sigma_n(e_1)\ \sigma_n(e_{a_1}))
(\sigma_n(e_2)\ \sigma_n(e_{a_2}))
\cdots
(\sigma_n(e_r)\ \sigma_n(e_{a_r})).
$$
Under this bijection, there is an edge from $\sigma_n(i)$
to $\sigma_n(j)$ in the graph representing $\sigma_n \phi^{-1}\sigma_n$
just when there is an edge from $i$ to $j$ in the graph representing
$\phi$.  Hence this bijection preserves treeness or unicyclicity of the
representation.  Setting $\phi=\pi$ or $\phi=\pi\mu$, this
proves (\ref{e1}) and (\ref{e2}).
Finally, we need to show (\ref{e3}).  If we apply (\ref{eqo}) to both
of its sides, this reduces to showing that 
$\Nb_k(\sigma_n \pi^{-1}\sigma_n)=\Nb_k(\pi)$ for all $k$.  
Let the reverse of a sequence $A$ be $\rev{A}$, and write
$\sigma_n(A)$ to denote the sequence derived by applying
$\sigma_n$ to each element of $A$.  Then, if
$\pi=(A)$, 
$\pi^{-1}$ will be $(\rev{A})$, and
$\sigma_n \pi^{-1} \sigma_n$ will be $(\sigma_n(\rev{A}))$,
so we wish to show that
$$
\sum_{(A_1\cdots A_k)=(A)} \prod_i \Ntz((A_i\  0)) =
\sum_{(A'_1\cdots A'_k)=(\sigma_n(\rev{A}))} \prod_i \Ntz((A'_i\  0)).
$$
There is a bijection between terms on the left- and right-hand sides
of this equation given by setting $A'_i=\sigma_n(\rev{A}_{k+1-i})$.
It will, therefore, do to show that 
\begin{eqnarray}
&&\Ntz((D\ 0))=\Ntz((\sigma_n(\rev{D})\  0)),
\\
\label{e4}
&&\mbox{for all sequences $D$ of distinct 
elements contained in $\{1,\ldots,n\}$.}\nonumber
\end{eqnarray}

Let the range of $D$ be $\{d_1<d_2<\cdots<d_q\}$.
$\Ntz((D\ 0))$ is the number of ways
of writing $(D\ 0)$ as an exchange permutation without an exchange $(0\ x)$,
i.e., the number of ways of representing $(D\ 0)$ as
\begin{equation}
\label{eqx}
(D\ 0)=(d_q\ a_q)\cdots (d_2\ a_2)(d_1\ a_1), \qquad
a_1, \ldots, a_q\in\{0,d_1,d_2,\ldots,d_q\}.
\end{equation}
Extend $\sigma_n$ to $\{0,1,\ldots,n\}$ by setting $\sigma_n(0)=0$.
If we invert both sides of (\ref{eqx}) and conjugate it by $\sigma_n$,
we get
$$
(\sigma_n(\rev{D})\ 0)=(\sigma_n(d_1)\ \sigma_n(a_1))
\cdots (\sigma_n(d_{n-1})\ \sigma_n(a_{n-1}))(\sigma_n(d_n)\ \sigma_n(a_n)).
$$
Since $\sigma_n$ is order-reversing on $\{d_1,\ldots,d_q\}$, this gives
a bijection between the representations of $(D\ 0)$ and those of
$(\sigma_n(\rev{D})\ 0)$.  This proves (\ref{e4}), and hence Theorem
\ref{resb}.
\qed\medskip

\section{The most likely permutation, for all $n$}\label{fan}

After the following two lemmas,
we will be able to answer Robbins and Bolker's question by determining,
for each $n$, which 
permutation on $\{1,\ldots,n\}$ has greatest multiplicity.  

Given a sequence $D$ of distinct positive integers, let 
$\Nb'_k((D))$ be $\Nb_k((D))$ with the restriction that, in the sets
of tree representations that this expression counts, we only count
tree representations containing a transposition $(d\ 0)$, where $d$
is the smallest element of $D$.
\begin{lemma}
\label{sublemma}
Let $k$ be given.
For a fixed-length sequence $D$ of distinct positive
integers, the quantities 
\begin{enumerate}
\item $\Nb'_k((D))$, 
\item $\Nb_k((D))-\Nb'_k((D))$, and
\item $\Nb_k((D))$ 
\end{enumerate}
are each separately maximized when
$D$ is a cyclic permutation of a decreasing sequence.
Moreover, for 1 and 3, if $k=1$, this is the only case where the
maximum occurs.
\end{lemma}
\begin{proof}
Without loss of generality, let $D$ end in its smallest element, which we
call $d$.  We claim that
\begin{equation}
\label{eqsl1}
\Nb_k((D))-\Nb'_k((D))=\sum_{D' D''=D} \Nb_k((D')) \Nt{d}((D''))
\end{equation}
and
\begin{equation}
\label{eqsl2}
{1\over k}\Nb'_k((D))=\sum_{J D_1\cdots D_k=D}
\Ntz((J\ 0)) \Ntz((D_1\ 0)) \cdots \Ntz((D_{k-1}\ 0)) \Nt{d}((D_k)),
\end{equation}
where in these sums $J$, but not $D'$, $D''$, or the $D_i$'s, may be empty.
The lemma will follow from these claims.

To see where $\Nb'_k((D))$ is maximized, use (\ref{eqsl2}).
If $D$ is decreasing, then $J\ 0$, $D_1\ 0$, \dots, $D_{k-1}\ 0$, and $D_k$
will also be decreasing, so by Lemma \ref{rb}, each term on the right-hand side
of (\ref{eqsl2}) will be maximized.  If $k=1$, then the term $\Nt{d}((D))$
will appear on the right-hand side of (\ref{eqsl2}).  By Lemma \ref{rb},
this will be maximum only when $D$ is a cyclic permutation of a decreasing 
subsequence, so (\ref{eqsl2}) is maximized only when claimed.

$\Nb_k((D))$ and $\Nb_k((D))-\Nb'_k((D))$ must 
be handled together by induction on the length of $D$.  For each possible
length of $D$, if $\Nb'_k((D))$ and $\Nb_k((D))-\Nb'_k((D))$ are maximized 
when $D$ is decreasing, $\Nb_k((D))$ obviously will be as well,
and, if $k=1$, this will be the only place the maximum occurs, since it is
the only place where $\Nb'_k((D))$ is maximized.
For $\Nb_k((D))-\Nb'_k((D))$, use
(\ref{eqsl1}).  The induction hypothesis and Lemma \ref{rb} will then show that
the right-hand side of (\ref{eqsl1}) is maximized when $D$ is decreasing.

We now prove our claims.

\noindent{\it Proof of (\ref{eqsl1})}.
To prove (\ref{eqsl1}), look at its left-hand side.  It counts $k$-tuples
$(T_1,\ldots,T_k)$ of trees rooted at 0, representing $(D_1\ 0)$,
\dots, $(D_k\ 0)$, respectively, where $(D_1\ \cdots\ D_k)=(D)$ and
$T_1$, \dots, $T_k$ do not contain the transposition $(d\ 0)$.  The
$T_i$'s must contain a unique transposition $(d\ e)$, where $e\ne d$;
if this transposition is in $T_j$, removing it from $T_j$ will leave us with
a tree $U$ rooted at $d$ and a tree $V$ rooted at 0 which jointly represent
$(D_j\ 0)(d\ e)$.  Then $d$ and 0 must appear in separate cycles in
$(D_j\ 0)(d\ e)$, so $d$ appears after $e$ in $D_j$.  
Recall that the sequence $D$ ends with the element $d$.
Let $D'$ be
the initial segment of $D$ ending in $e$, and let $D''$ be the remainder
of $D$, so $D=D'\ D''$.  
Since $d$ appears after $e$ in $D_j$, it follows that $D''$ is a
subsequence of $D_j$.
Delete $D''$ from $D_j$ to give $D'_j$.  Now
$$(D_j\ 0)(d\ e)=(D'_j\ 0)(D''),$$
so $U$ represents $(D'')$ and $V$ represents $(D'_j\ 0)$.  In fact,
$T_1$, \dots, $T_{j-1}$, $V$, $T_{j+1}$, \dots, $T_k$ represent
$(D_1\ 0)$, \dots, $(D_{j-1}\ 0)$, $(D'_j\ 0)$, $(D_{j+1}\ 0)$,
\dots, $(D_k\ 0)$, respectively, so $(T_1,\ldots,T_{j-1},V,T_{j+1},\ldots,
T_k)$ is a representation of the form counted by $\Nb_k((D'))$.  We have
already remarked that $U$ is a tree rooted at $d$ representing $D''$.
This gives a map from the set counted by the left-hand side of (\ref{eqsl1})
to the set counted by the right-hand side of (\ref{eqsl1}).  This map
is easily seen to be invertible, proving (\ref{eqsl1}).
{\qed \medskip}

\noindent{\it Proof of (\ref{eqsl2})}.
The left-hand side of (\ref{eqsl2}) counts the same thing as the left-hand
side of (\ref{eqsl1}), except that the transposition $(d\ 0)$ is required
instead of forbidden.  Also, we assume that $D_k$ contains $d$; this gives
the factor of ${1\over k}$.  Let $D_k=E\ J$, where $E$ ends with $d$.  If
the transposition $(d\ 0)$ is deleted from $T_k$, we will be left with 
a tree $U$ rooted at $d$, representing $(E)$, and a tree $V$ rooted at 0,
representing $(J\ 0)$.  However, $J\ D_1\ \cdots\ D_{k-1} E=D$, so the
$(k+1)$-tuple of trees $(V,T_1,\ldots,T_{k-1},U)$ is of the form counted
by the right-hand side of (\ref{eqsl2}).  This gives us a map from the
set counted by the left-hand side of (\ref{eqsl2}) to the set counted by
the right-hand side of (\ref{eqsl2}).  This map is invertible, so (\ref{eqsl2})
holds.
\end{proof}

\begin{lemma}
\label{claim2}
Among all products of two disjoint cycles of fixed length, $m$ and $p$,
acting on a set of cardinality $m+p$, the permutations
\begin{equation}
\label{eqxxx}
(c_1\ c_2\ \cdots\ c_m)(c_{m+1}\ c_{m+2}\ \cdots\ c_{m+p})
\end{equation}
and
\begin{equation}
\label{eqxyy}
(c_1\ c_2\ \cdots\ c_p)(c_{p+1}\ c_{p+2}\ \cdots\ c_{m+p}),
\end{equation}
where $$c_1>c_2>\cdots>c_{m+p},$$
have the highest unicyclic multiplicity.  
This  multiplicity is
\begin{equation}\label{2cycles}
\Nu((m+p\ \cdots\ p+2\ p+1)(p\ \cdots\ 2\ 1))
=\frac{\bbb{m} \bbb{p} (m+p+4mp)-(m+p){{2m+2p}\choose{m+p}}}
{2(m+p)(m+p+1)}.
\end{equation}
%
% Do we wish to show WHY G_{mp} equals this?
%
All other such products have smaller
unicyclic multiplicity.
\end{lemma}
\noindent{{\it Remark}.  The proof of the exact formula is left to the reader as 
a tedious exercise based on (\ref{eqo}).}\hfill\\

\begin{proof}
Let $A=a_1\ \cdots\ a_m$ and $B=b_1\ \cdots\ b_p$ be sequences of distinct
positive integers with disjoint ranges, and let $\pi=(A)(B)$.
We wish to show that, for fixed $m$ and $p$,
$\Nu(\pi)$ is maximized iff we can write
\begin{eqnarray}
\label{cl2mx}
&&(A)=(c_1\ \cdots\ c_m)\ \ \mbox{ and }\ \ 
(B)=(d_1\ \cdots\ d_p), \ \ \mbox{where}\\
&&\mbox{$c_1>\cdots>c_m$, $d_1>\cdots>d_p$, and either $c_m>d_1$ or $d_p>c_1$.}
\nonumber
\end{eqnarray}
We will induce on $m+p$.  Without loss of generality, let 1 be the
smallest element in the domain of $\pi$.

Exchange $m$ and $p$, if necessary, to put
1 in $B$, so that we can set $b_p=1$.  In the case (\ref{cl2mx}), then,
we must have $d_p=1$, so $c_1>\cdots>c_m>d_1>\cdots>d_p$.
Let $H$ be a unicyclic representation of $\pi$, with ring $K$.
$H$ will contain some transposition $(1\ j)$, $j\ne 1$.  If 1 is not in $K$, 
then deleting the edge from 1 to $j$ from $K$ will leave a graph with a tree 
component $T$ rooted at 1 and a unicyclic component $H'$.  This graph
will represent $\pi(1\ j)$.  Considering the postman for a moment, we
see that in $H$,
after stopping for the night at $j$, he must proceed immediately to 1
and traverse $T$.  Hence $j$ and 1 are in the same cycle of $\pi$, and
$$\pi(1\ j)=(a_1\ \cdots\ a_m)(b_1\ \cdots\ b_q=j)(b_{q+1}\ \cdots\ 
b_p=1),$$
for some $1\le q\le p-1$.  Evidently, $T$ must represent $(b_{q+1}\ \cdots\ 
b_p=1)$, and $H'$ must represent 
$(a_1\ \cdots\ a_m)(b_1\ \cdots\ b_q)$.  Conversely, given a tree rooted
at 1 representing $(b_{q+1}\ \cdots\ b_p)$, and a unicycle representing
$(a_1\ \cdots\ a_m)(b_1\ \cdots\ b_q)$, we can insert an edge from 1 to $b_q$
to give a unicyclic representation of $\pi$.  Hence
\begin{equation}
\label{cl2ref}
Q(\pi)=\sum_{1\le q\le p-1} \Nt 1((b_{q+1}\ \cdots\ b_p=1))
\Nu((a_1\ \cdots\ a_m)(b_1\ \cdots\ b_q)),
\end{equation}
where $Q(\pi)$ is the portion
of $\Nu(\pi)$ coming from graphs with 1 not in $K$.  
By Lemma \ref{rb},
the $\Nt 1$ term will be maximized just when
$b_{q+1}>\cdots>b_p=1$, which will certainly happen in the case (\ref{cl2mx}),
and in this case, by the induction
hypothesis, the $\Nu$ term will be maximized as well.
Hence $Q(\pi)$ is never any bigger than it is in the case (\ref{cl2mx}).

Let $Q'(\pi)=\Nu(\pi)-Q(\pi)$ count the remaining unicyclic representations, 
that is, those where 1 is in $K$.
Let there be edges in $K$ from $z$ to 1 and from 1 to $w$.
Since $z>1$ and $w>1$, 1 is lower and $w$ is upper, so $(A)$ must be the
upper cycle and $(B)$ the lower.  Also, using the notation
of the proof of Proposition \ref{claim1}, 
the last element of some 
$E_j$ is 1; that is, $1=\bm(U_j)$ for some $j$. 
Conversely, if $1=\bm(U_j)$
for some $j$, 1 will certainly appear in $K$.  Therefore, we can think
of $Q'(\pi)$ as the right-hand side of (\ref{eqm}), if we apply the
following constraint: in the sets of tree representations that this
right-hand side counts, delete all tree representations not containing
a tree $U_j$ with $1=\bm(U_j)$.

In the case (\ref{cl2mx}), the inequalities between $\bp$, $\gp$,
$\bm$, and $\gm$ in (\ref{eqm}) are vacuous, so we may remove them,
and think of $Q'(\pi)$ as the right-hand side of (\ref{eqo}), subject
to our constraint.  
In all other cases, relaxing these inequalities 
does not decrease the right-hand side of (\ref{eqm}), so it will
do to prove that the right-hand side of (\ref{eqo}), subject to
our constraint, is never bigger than in the case (\ref{cl2mx}).  

Our constraint is equivalent to demanding that a transposition $(1\ 0)$
occur somewhere in the tree representations of the $B_j$'s, so the
right-hand side of (\ref{eqo}) with our constraint applied is
\begin{equation}
\label{eqzz}
\sum_{k\ge 1} 
\frac{1}{k}
\Nb_k((A)) \Nb'_k((B)).
\end{equation}
It follows from Lemma \ref{sublemma} that (\ref{eqzz}) is maximized in
the case (\ref{cl2mx}).  We still need to show exactly when the maximum occurs.
For all $D$, $\Nb_1((D))\ge\Ntz((D\ 0))$, and $\Ntz((D\ 0))$ is positive
by \cite[Theorem 4]{ref4}.
Define $\Ntzprime((D\ 0))$ to be the number of tree representations of $(D\ 0)$ 
with 0 as the root and a transposition $(d\ 0)$ present, for $d$ the
smallest element of $D$.  Then, for any $D$, let $D'=D''\ d$ be $D$ cyclically
rotated so that its smallest element, $d$, is at the end.  In this case,
$\Nb'_1((D))\ge \Ntzprime((D''\ d\ 0))$.
As in the proof of Lemma \ref{sublemma}, if a 
tree representation of $(D''\ d\ 0)$ rooted at 0 contains 
the transposition $(d\ 0)$, we can delete it to get a tree representation
of $(D''\ d)$ rooted at $d$.  This transformation is reversible, so
$$\Ntzprime((D''\ d\ 0))=\Nt d((D''\ d)),$$
and by \cite[Theorem 4]{ref4}, 
this is positive.  Hence $\Nb_1((D))$ and $\Nb'_1((D))$
are always positive. 
Also, by Lemma \ref{sublemma},
we know that the maximum of $\Nb'_1((D))$ and $\Nb_1((D))$ 
can only occur 
when $D$ is a cyclic permutation of a decreasing sequence.
Looking at the $k=1$ term in (\ref{eqzz}), then, we see that for the maximum
to occur, we must be able to write
\begin{eqnarray*}
&&(A)=(c_1\ \cdots\ c_m)\ \ \mbox{ and }
\ \ (B)=(d_1\ \cdots\ d_p), \ \ \mbox{where}\\
&&\mbox{$c_1>\cdots>c_m$ and $d_1>\cdots>d_p=1$.}
\end{eqnarray*}
If $p=1$, we are done.  Otherwise, write $B=B'\ 1$, and 
look at the $q=p-1$ term in (\ref{cl2ref}).  By the induction hypothesis,
if we are at the maximum, (\ref{cl2mx}) must be satisfied for $A$ and $B'$.
Therefore, either $c_m>d_1$, and we are done, or $d_{p-1}>c_1$.

If $d_{p-1}>c_1$, look at the representations of $\pi$ with 1 in $K$.
Each $\gm(B_j)$ must be exceeded by a $\bp(T_i)$; therefore, 1 is the only
possible value for any $\gm(B_j)$, so $k=1$ in (\ref{eqm}) and $B_1=1\ B'$.
Ignoring the other conditions on $\bm$, $\gm$, $\bp$, and $\gp$ that occur 
in (\ref{eqm}), we see that
$$
Q'((A)(B))\le \Nb_1((A)) \Ntzprime((1\ B'\ 0)).
$$

Let $A'$ denote a shifted version of $A$ with every element greater than
every element of $B$.  Then, as we know,
$$
Q'((A')(B))\ge\Nb_1((A'))\Nb'_1((B))=\Nb_1((A))\Nb'_1((B)),
$$
and since we know that $\Nb_1((A))>0$, it will do to show that
$\Nb'_1((B))>\Ntzprime((1\ B'\ 0)).$  Now,
$\Nb'_1((B))$ is a sum of $\Ntzprime((B_1\ 0))$ over cyclic
permutations $B_1$ of $B$.
Since $p>1$, $1\ B'\ne B'\ 1$, so it will do to show that 
$\Ntzprime((B'\ 1\ 0))>0$; but we know this from above.
This shows that the case $d_{p-1}>c_1$ cannot occur, completing the proof.
\end{proof}

Let $G_{mp}$ denote the right-hand side of (\ref{2cycles}).
\begin{corollary}
\label{corrdx}
Fix a conjugacy class of $S_n$, and let it consist of all permutations
whose representations as products of disjoint cycles consist of an $L_1$-cycle,
an $L_2$-cycle, \dots, and an $L_q$-cycle, where $L_1+\cdots+L_q=n$.
Then the permutations of maximum
multiplicity in this class are just the permutations where each cycle of 
length $L$ is of the form $(\Sigma+L\ \cdots\ \Sigma+1)$, for some $\Sigma$.
Furthermore, the multiplicity of the permutations of maximum multiplicity is
\begin{equation}
\label{eqmaxx} 
\sum_{\shortstack{\scriptsize $\chi$ an involution of\\
$\scriptstyle \{1,\ldots,q\}$}}\ 
\prod_{\chi(i)=i} C_{L_i}
\prod_{\chi(i)=j,\ i<j} G_{L_i L_j}.
\end{equation}
\end{corollary}
\begin{proof} 
Let $\pi$ be an arbitrary permutation of $\{1,\ldots,n\}$, and let
$\pi=\pi_1\cdots\pi_q$ be its decomposition into disjoint cycles, where
$\pi_1$ is an $L_1$-cycle, $\pi_2$ is an $L_2$-cycle, \dots,
and $\pi_q$ is an $L_q$-cycle.  By Lemma \ref{key},
$$
\N(\pi)=\sum_{\shortstack{\scriptsize $\chi$ an involution of\\
$\scriptstyle \{1,\ldots,q\}$}}\ 
\prod_{\chi(i)=i} \Ntree(\pi_i)
\prod_{\chi(i)=j,\ i<j} \Nu(\pi_i \pi_j).
$$
Using the bounds in Lemmas \ref{rb} and \ref{claim2} then shows that
$\N(\pi)$ is no bigger than (\ref{eqmaxx}) and that the equality condition
is as claimed.
\end{proof}

Let \emph{cycle left} be the permutation 
$(n\ n-1\ \cdots\ 2\ 1),$ and 
let {\em double cycle left} be the permutation 
$$(n\ n-1\ \cdots\ m+1)
(m\ \cdots\ 2\ 1), 
$$
where $m$ is either $\floor{n/2}$ or $\ceil{n/2}$.
Note that the multiplicity of double cycle left is given by Lemma~\ref{claim2}.

\begin{theorem}
For fixed $n$, the permutations with maximum multiplicity are:
for $n=1$, the identity; 
for $n=2$, cycle left or identity;
for $n=3$, cycle left or double cycle left;
for $4\le n\le 17$, double cycle left;
for $18\le n$, the identity.
\end{theorem}
\begin{proof}
For $n\ge 29$, we already know this from Proposition \ref{idbd1},
so fix some $n\le 28$.  Let $(L_1,L_2,\ldots,L_q)$ be a partition of $n$
into positive integers.  We can compute (\ref{eqmaxx}) for each such 
partition.  The maximum value
is found to be achieved for $q=n$ and $L_1=\cdots=L_n=1$ for $n=1$ and 
for $18\le n\le 28$;
for $q=2$ and $\{L_1,L_2\}=\{\floor{n/2},\ceil{n/2}\}$ for $4\le n\le 17$;
and for either $q=1$ and $L_1=n$ or $q=2$ and 
$\{L_1,L_2\}=\{1,n-1\}$ for $n=2$ and $n=3$.  Applying Corollary
\ref{corrdx} completes the proof.
\end{proof}

\section{The limiting distribution of the number of fixed points}\label{fix}

Schmidt and Simion \cite{ref6} consider the exchange permutation
(\ref{eone}) where $a_1$, \dots,
$a_n$ are picked independently at random from the uniform distribution
on $\{1,\ldots,n\}$.  They ask what the limiting probability $p_k$ of
$\pi$ having $k$ fixed points is, as $n$ becomes large.

\begin{theorem}
\label{ss1}
If $a_1$, \dots, $a_n$ are picked independently uniformly at random,
the limiting probabilities $p_k$ of $\pi$
given by (\ref{eone}) having exactly $k$ fixed points
as $n\rightarrow\infty$ are given by
$$
\sum_{k\ge 0} p_k u^k=q(u),
$$
where
$$q(u)=\exp({1\over 2}(u-1)(u e^{-2}-6e^{-1}+4-e^{-2})).$$
\end{theorem}
First 
we recall the Lagrange Inversion Theorem (see \cite{ref2}), 
and give three examples
that will be used in the proof.
For a formal power series $\phi(z)=\phi_0 +\phi_1z +\phi_2z^2 + \cdots\ $, let
$[z^n]\phi$ denote $\phi_n$, the coefficient of $z^n.$
\begin{theorem}(Lagrange Inversion Theorem, or LIT).
Let $\phi,\psi,$ and $A$ be formal power series with complex coefficients,
with  $\phi_0\neq0.$ If $\psi(z)=z\phi(\psi(z))$, then
$$[z^n]A(\psi(z))=\frac1n[z^{n-1}]\phi^n(z)A'(z)$$
for any positive integer $n$.
\end{theorem}

Here are some applications of LIT that will be used to prove Theorem~\ref{ss1}.
\begin{example} The number $\mathrm{Tr}_n$ of 
rooted trees on $n$ labeled vertices.
Let $t(z)$ be the exponential generating function for $\mathrm{Tr}_n$,
$t(z) = \sum_{n>0}\mathrm{Tr}_n z^n/n!.$ One checks  that 
$t=ze^t$. Apply LIT with $\psi(z)=t(z), \phi(z) = e^z,$ 
(note that the constant term of $e^z$ is nonzero) and 
$A(z)=z$ to get $\mathrm{Tr}_n=n^{n-1}.$
\end{example}
\begin{example} The expansion of $z/t(z).$
Apply LIT with $\psi(z)=t(z), \phi(z) = e^z,$ and 
$A(z)=e^{-z}$ to get  $$\frac{z}{t(z)}=1-\sum_{n>0}\frac{(n-1)^{n-1}}{n!}z^n.$$
\end{example}
\begin{example} The expansion of $1/(1-t(z)).$
Apply LIT with $\psi(z)=t(z), \phi(z) = e^z,$ and 
$A(z)=\frac{1}{1-z}$ to get  
$$\frac{1}{1-t(z)}=\sum_{n\ge0}\frac{n^{n}}{n!}z^n.$$
\end{example}
\noindent{\it Proof of Theorem \ref{ss1}}.
Let $T_{nk}$ be the number of $(a_1,\ldots,a_n)$'s inducing a connected digraph
and a permutation with exactly $k$ fixed points.
Using \cite[Lemma 3.5]{ref6} 
gives $\sum_{k>0} T_{1k}=\sum_{k>0} T_{2k}=1$, and
\mbox{$\sum_{k>0} T_{nk}=2(n-1)^{n-1}$} for $n\ge 3$.
We claim that $T_{11}=T_{22}=1$, $T_{n1}=2(n-1)^{n-1}$ for $n\ge 3$, and
$T_{nk}$ is zero for all other positive values of $k$.
This is because  a connected digraph has exactly zero or one fixed points, 
unless $n=2$, when it must have zero or two.  
For  $t(z)$ as in  Example~1, 
we can then  write
$$
T(z,u)=\sum_{n>0,k>0} T_{nk} {z^n u^k\over n!} = 
zu+{1\over 2} z^2u^2 + u(2-{2z\over t(z)}-2z-z^2),
$$
by Example~2.
The exponential generating function for the number of ways to pick 
$(a_1,\ldots,a_n)$ overall is
$$
U(z)=\sum_{n\ge 0} {n^n\over n!}z^n.
$$
Observe that $\log U(z)$ is the 
exponential generating function for 
the number of connected digraphs on $n$ vertices.
The identity $f(z,u)=f(z,1)$, true for any function $f$ which 
does not depend on the variable $u$, holds for 
$$f(z,u)=-T(z,u)+\sum_{n>0,k\ge 0} T_{nk}{z^n u^k\over n!}.$$
However, when $u$ is set equal to 1, this equals $-T(z,1)+\log U(z)$,
so the exponential generating function
for components where $u$ marks the number of fixed points is given by
$$
\sum_{n>0,k\ge 0} T_{nk}{z^n u^k\over n!}=T(z,u)-T(z,1)+\log U(z).
$$

Let $p_{nk}$ be the probability of exactly $k$ fixed points occurring in a permutation
on $\{1,\ldots,n\}$, and let $q_n(u)=\sum_{k\ge 0} p_{nk} u^k$.  Then
$$
q_n(u)={[z^n] \exp(T(z,u)-T(z,1)+\log U(z))\over [z^n] U(z)}
={[z^n] U(z)\exp(T(z,u)-T(z,1))\over [z^n] U(z)}.
$$
By Example 3, $U(z)=1/(1-t(z))$.  
Now, since both numerator and denominator of $q_n(u)$
are analytic in the complex plane 
slit along $[e^{-1},\infty)$, we can use singularity analysis \cite{ref2a}
and the 
expansion \cite{ref1} at $z=e^{-1}$ 
$$
t(z)=1-\sqrt{2(1-ez)}-{1\over3}(1-ez)+O((1-ez)^{3/2})
$$
to find that $q_n(u)\rightarrow q(u)$
pointwise in $u$
as $n\rightarrow\infty$.  By \cite[Theorem 9.3]{ref3}, the limiting probability
$p_k$ of exactly $k$ fixed points as $n\rightarrow\infty$ is then given by
$\sum_{k\ge 0} p_k u^k=q(u)$.
\qed
\medskip

Schmidt and Simion \cite{ref6} also ask what the limiting distribution 
as $n\to\infty$ of the
number of fixed points
is if $a_1$, \ldots, $a_n$ are  chosen to be a random permutation
of $\{1,\ldots,n\}$.

\begin{theorem}
If $a_1$, \dots, $a_n$ are picked to be a random permutation of 
$\{1,\ldots, n\}$,
the limiting probabilities $\bar p_k$ of $\pi$
given by (\ref{eone}) having exactly $k$ fixed points
as $n\rightarrow\infty$ are given by
$$
\sum_{k\ge 0} \bar p_k u^k=\bar q(u),
$$
where
$$
\bar q(u)=\exp({1\over 2}(u-1)(u+4e-7)).
$$
\end{theorem}
\begin{proof}
Define $\bar T_{nk}$, $\bar T(z,u)$, etc. analogously
to the definitions of $T_{nk}$, $T(z,u)$, etc. in the proof of
Theorem \ref{ss1}.  Then
$\bar T_{nk}=0$ for $k>0$ except that
$\bar T_{11}=1$, $\bar T_{22}=1$, and, for $n\ge 3$, $\bar T_{n1}=2$
(cf. \cite[Theorem 3.13]{ref6}); therefore, 
$$\bar T(z,u)=\sum_{n>0,k>0} \bar T_{nk} {z^n u^k\over k!} = 
zu+{1\over 2} z^2u^2 + 2u(e^z-1-z-{1\over 2}z^2)$$
and
$$
\bar U(z)=\sum_{n\ge 0} {n!\over n!}z^n={1\over 1-z}.
$$
As in Theorem \ref{ss1},
$$
\bar q_n(u)={[z^n] \bar U(z)\exp(\bar T(z,u)-\bar T(z,1))\over [z^n] \bar U(z)}.
$$
Both numerator and denominator are analytic except for a
pole at $z=1$, so using singularity analysis \cite{ref2a} there we get
$\bar q_n(u)\rightarrow \bar q(u)$ pointwise in $u$
as $n\rightarrow\infty$, and by \cite[Theorem 9.3]{ref3}, the 
generating function of
the limiting probabilities is $\bar q(u)$.
\end{proof}

\enlargethispage*{1000pt}
\pagebreak

\enlargethispage*{1000pt}
\pagebreak

\section*{Frequently used notations}

\medskip
\begin{tabular}{ll}
\vspace{10pt}
Notation & Definition\\
$\Ntree(\pi)$&the number of representations of $\pi$ as a tree\\
$\Nu(\pi)$ &the number of  representations of $\pi$ as a unicycle\\
$\Nt x(\pi)$ & the number of tree representations of $\pi$\\ 
& where  $a_x=x$\\
$\Ntz(\pi)$ & $\Nt 0(\pi)$\\
$\bp$ (resp. $\bm$)  & the largest (resp. smallest) neighbor of $0$\\
$\gp$ (resp. $\gm$) &the last (resp. first) element of a sequence \\
$\Ntz_{\bp>x}(\pi)$ & the number of representations
of $\pi$ as a \\
&tree $T$ with root 0 and $\bp(T)>x$\\
$\Ntz_{\bm<x}(\pi)$ & the number of representations
of $\pi$ as a \\
& tree $T$ with root 0 and $\bm(T)<x$\\
$\Nuu((A),(B))$ & the number of  unicyclic representations\\
&of the permutation $(A)(B)$  for which $(A)$\\ 
& is the upper cycle\\
$\Nb_k((D))$ & $\sum_{(D_1\cdots D_k)= (D)} \prod_i \Ntz((D_i\  0))$\\
$\pi\BB\rho$ & the concatenation of permutations $\pi$ and $\rho$ \\
$\sigma_n$ & the permutation sending $i$ to $n+1-i$ \\
$\flip{\pi}$ & $\sigma_n \pi^{-1} \sigma_n$, if $\pi$ acts on $\{1,\ldots,n\}$\\
$\rev{A}$ &  the reverse of a sequence $A$\\
$Q(\pi)$ &the portion of  $\Nu(\pi)$ coming from graphs\\
&with 1 not in the ring \\
$Q'(\pi)$ &the portion of  $\Nu(\pi)$ coming from graphs\\
&with 1 in the ring \\
$\Nb'_k((D))$ &the portion of $\Nb_k((D))$ with $a_d=0$,\\
&where $d$ is the smallest element of $D$\\
$\Ntzprime((D\ 0))$ &the number of tree representations of $(D\ 0)$\\ 
&with $a_0 = 0$ and  $a_d=0$, where $d$ is the \\
&smallest element of $D$
\end{tabular}

\vspace{.75in}

\noindent Center for Communications Research\\
4320 Westerra Court\\
San Diego, CA 92121\\
{\tt dgoldste@ccrwest.org}\\
{\tt dmoews@ccrwest.org}

\end{document}